\documentclass[12pt,a4paper]{article}

\usepackage[english]{babel}
\usepackage{graphicx}
\usepackage[centertags]{amsmath}
\usepackage{amsfonts}
\usepackage{amssymb,amsbsy}
\usepackage{amsthm}
\usepackage{newlfont}

\newtheorem{thm}{Theorem}
\newtheorem{cor}[thm]{Corollary}

\newtheorem{prop}[thm]{Proposition}
\newtheorem{rem}{Remark}

\newcommand{\be}{\beta}
\newcommand{\la}{\lambda}
\newcommand{\wt}{\widehat}

\newcommand{\oln}{\overline}
\newcommand{\uln}{\underline}

\newcommand{\mcm}{\mathcal M}
\newcommand{\mbm}{\mathbf m}
\newcommand{\mby}{\mathbf y}
\newcommand{\mbx}{\mathbf x}
\newcommand{\mbz}{\mathbf z}
\newcommand{\mbb}{\mathbf b}

\newcommand{\mbq}{\mathbf q}
\newcommand{\mbf}{\mathbf f}

\newcommand{\nn}{\noindent}
\def\ind#1{{\mbox{\fontsize{9}{10.8}\selectfont$\relax#1 $} }}

\begin{document}

\title{Fast algorithms of Bayesian Segmentation of Images}

\author{Boris A. Zalesky\\
Institute of Engineering Cybernetics NAN, Minsk, Belarus
\footnote{address: 220012, Surganov str. 6, Minsk,
Belarus}\\
zalesky@mpen.bas-net.by
}



\maketitle

\begin{abstract}
The network flow optimization approach is offered for
Bayesian segmentation of  gray-scale and color images.
It is supposed image pixels are characterized by a
feature function taking finite number of arbitrary rational values
(it can be either intensity values or other characteristics of images).
The clusters of homogeneous pixels are described by
labels with values in another set of rational numbers.
They are assumed to be dependent and distributed
according to either the exponential or the Gaussian Gibbs law.
Instead traditionally used local neighborhoods of nearest pixels
the completely connected graph of dependence of all pixels
is employed for the Gibbs prior distributions.

The methods developed reduce the problem of segmentation
to the problem of determination of the minimum cut of an appropriate network.
\end{abstract}

\vspace{3cm}
\nn

\text{}\\
\nn

\newpage

\section{Introduction}

The segmentation of images and restoration of degraded images are branches of
image processing that are now extensively studied for
their evident practical importance as well as
theoretical interest. There are many approaches to
solution of these problems.
We consider here methods of the Bayesian segmentation of images,
which generalize methods of Bayesian image estimation
\cite{GG84,Gi95,Z00}.
Because of  use of a prior information the methods of Bayesian image
estimation
would be methods of the first choice for many
practical problems but unfortunately
they  are often difficult to compute.
Until recently
approximations of the Gibbs estimators  have been usually available
rather then their exact values.
For last decade the significant progress in the Gibbs estimation
have been achieved. In particular, the methods of discrete
optimization for the high-dimensional Gibbs
estimation and segmentation have been obtained \cite{BVZ01,KZ01,PRIP01,DAN01,AMS01}.
The methods developed allowed not only an
efficient evaluation of the Gibbs
estimates and but also their
fast exact determination.

In this paper the  problem of image segmentation
is considered as a problem of cluster-analysis  for the case of
Gibbs prior distributions of clusters.
There are
enormous number of papers devoted to the problems of
Bayesian cluster-analysis (see \cite{ABEM89,McL92,MM99,JH98}).
In some of them methods of classification with dependent
clusters are presented \cite{MM99,JH98}.
We consider the problem of classification of observations
$y_{1},y_{2},$ $\ldots,y_{n}$ (for instance, image intensities or
texture characteristics of images etc.)
with the feature function $f(y_i)$ that takes finite
number of rational values. The clusters
are identified by appropriate rational labels
(note at once that in the presented below models the case of
rationaly valued feature functions and cluster labels is reduced
to the case of functions and labels taking only integer values).
The clusters are supposed  dependent and distributed according to the Gibbs
field (such type models occur frequently in image processing).
The labeling Gibbs field is specified on the directed fully
connected graph of all pixels
(usually only graphs of connected nearest neighbors pixels were
considered. Often it was finite d-dimensional regular lattices).

The discrete
optimization methods that enable efficient determination of an
exact solution of the segmentation problem are described. In
the offered methods the problem of identification of the Gibbs
classifiers  is reduced, first, to the integer optimization
problem and then to the problem of findind of the minimum cuts of
special networks. The minimum cuts of the networks can be found by
fast methods \cite{AMS01,CG94} that  take into consideration the
specific character of the networks and allow
execution in concurrent mode.

\section{Description of Models}

It can be easily seen that in   models presented below the case of
rationally valued feature functions and cluster labels is reduced
to the case of integer classification.
So, let  an image
$\mby=(y_{1},y_{2},\ldots ,y_{n})$ be with the
feature function $f$
that takes finite number integer values in
$Z_L=\{0,1,\ldots,L\}$ and let $f_i=f(y_i)$. Let
$\mathcal{M}=\{m(0),m(1),\ldots,m(k)\}$, $(0\le m(0)< m(1)<
\ldots< m(k)\le L)$ be a set of allowable cluster labels. Suppose
the image $\mby$ is partitioned into $k+1$
clusters (the number $k+1$ instead of $k$ is taken only
to simplify formulas)
and  each cluster is specified by integer number $j,\
(0\le j\le k)$, as well as by an appropriate fixed integer label
$m_j\in \mathcal{M}$.

The feature vector $\mbf=(f_{1},f_{2},\ldots ,f_{n})$
of image $\mby$ is considered as
a random variable with  the non-identical exponential
\begin{equation*}\label{e:p_e}
p_{i,exp}^{\, m}(f_i)=c_i\cdot \exp\left\{-\la_i|f_i-m_i| \right\}
\end{equation*}
or  with the non-identical Gaussian
\begin{equation*}\label{e:p_g}
p_{i,gaus}^{\,
m}(f_i)=c_i\cdot\exp\left\{-\la_i(f_i-m_i)^2\right\}
\end{equation*}
Gibbs distribution  $(\la_i \ge 0)$.

The labels $m_i$ of clusters are supposed to
be dependent random variables distributed according to the Gibbs
field. The segmented image
$\mbx=\{x_{1},x_{2},\ldots,x_{n}\}$
is specified on the fully connected directed  graph
$G=(V,E)$ with the set of vertices $V=\{0,1,\ldots,n\}$ and the
set of directed arcs $E=\{(i,j)\ |\ i,j\in V\}$. It takes values
in the space of labels $\mathcal{M}^V$, i.e. $x_i\in
\mathcal{M},\ i\in V$, and  is  either  the exponential or the
gaussian form
\begin{equation*}
\begin{split}
&{p}_{exp}(\mbm)=c\cdot\exp\{-\sum_{(i,j)\in E}\beta_{i,j}|m_i-m_j|\},\\
&{p}_{gaus}(\mbm)=c\cdot\exp\{-\sum_{(i,j)\in
E}\beta_{i,j}(m_i-m_j)^2\},
\end{split}
\end{equation*}
where $c$ is the norming quantities (here and below different
constant will be denoted by the same letter), the vector
$\mbm=\{m_{1},m_{2},\ldots ,m_{n}\}$ and parameters
$\beta_{i,j}\ge 0$.
\begin{rem} (i) The Gibbs models considered are extentions of the Ising model.

\noindent
(ii) Note, that for some image processing problems the
prior distribution ${p}_{exp}(x)$ is more preferable than
Gaussian Gibbs one because of smaller blurring effect. Moreover,
its identification turns out far more computationally efficient.
\end{rem}

For the model with the
exponential conditional distribution $p_{i,exp}^{\, m}(f_i)$ and
the exponential prior $p_{exp}(\mbm)$ the Gibbs classifier is equal to
\begin{equation}\label{e:m_e}
\widehat{\mbm}_{exp}=\arg\min_{m\in \mathcal{M}^V}
\left\{\la_i\sum_{i\in V}|f_i-m_i|+
\sum_{(i,j)\in E}\be_{i,j}|m_i-m_j|\right\},
\end{equation}
and respectively, for the model with the Gaussian conditional
distribution $p_{i,gaus}^{\, m}(f_i)$ and the Gaussian  prior
$p_{gaus}(m)$ the Gibbs classifier is of the form
\begin{equation}\label{e:m_g}
\widehat{\mbm}_{gaus}=\arg\min_{m\in \mathcal{M}^V}
\left\{\sum_{i\in V}\la_i(f_i-m_i)^2+ \sum_{(i,j)\in E}
\be_{i,j}(m_i-m_j)^2\right\}.
\end{equation}
In spite of clear posing the problem of finding exact values of
the Gibbs classifiers $\widehat{\mbm}_{exp}$ and $\widehat{\mbm}_{gaus}$ for
large samples is rather complicated. So, for instance,  computer
images, which  are frequent objects of classification, usually
consist up to $2^{18}$ variables or more. Nevertheless, it turned
out possible to compute efficiently (during polynomial time on
the sample size $n$ and number of clusters $k$)  both of these
classifiers.  Theoretically, run-time for the first classifier
does not exceed $ckn^3$ and for the second one is less than
$c(kn)^3$. For many applied problems real run-time was even of order
$O(nk)$. The efficient computation of the classifiers for the mixed model
with the exponential conditional distributions $p_{i,exp}^{\,
m}(f_i)$ and the Gaussian prior $p_{gaus}(\mbm)$ as well as for
the mixed model with the Gaussian conditional distributions
$p_{i,gaus}^{\, m}(f_i)$ and the exponential prior
$p_{exp}(\mbm)$ is also available.

\section{Computation of the Optimal Classifiers}

In the case of two classes ($k=1$) the Gibbs classifier $\wt{\mbm}_{exp}$ and
$\wt{\mbm}_{gaus}$ coincide (since any Boolean variable $b$
satisfies the identity  $b^2=|b|$). They can be evaluated by the
network flow optimization methods \cite{GPS89,PR75}. Moreover,
in 1989  Greig, Porteous and Seheult \cite{GPS89} developed a heuristic
network flow algorithm that is especially efficient for
estimating the Boolean Gibbs estimator. These authors   posed
also the problem for the case more than 2 clusters. Recently we have
described  the multiresolution network flow minimum cut
algorithm \cite{AMS01} that allows exact computation of Boolean classifiers as
well as developed algorithms of computation of the mentioned Gibbs classifiers in
general case. It turned out the network flow optimization methods
can be used to identify $\wt{\mbm}_{exp}$ and $\wt{\mbm}_{gaus}$
even when $k>1$.
\subsection*{Identification of $\wt{\mbm}_{exp}$} Denote the function
to be minimized by
$$
U_1(\mbm)=\sum_{i\in V}\la_i|f_i-m_i|+ \sum_{(i,j)\in E}
\be_{i,j}|m_i-m_j|, \quad (\mbm\in \mcm^V).
$$
The idea of the method is to represent the vector of  labels of
clusters $\mbm$ by the integer valued linear combination of Boolean
vectors and then reduce the problem of integer minimization of
the function $U_1(\mbm)$ to the problem of Boolean minimization.
The problem of Boolean minimization can be solved by the network flow
optimization methods.

Let for  arbitrary integers $\mu$ and $\nu$ the indicator function
$1_{(\mu\ge \nu)}$ be equal to 1 if $\mu\ge \nu$ and be equal to
0 otherwise. For any $\mu\in \mcm$ and Boolean variables
$x(l)=1_{\big(\mu\ge m(l)\big)}$
such that
$ x(1)\ge x(2)\ge\ldots\ge x(k)$
the identity
\begin{equation}\label{e:mu}
\mu=m(0)+\sum_{l=1}^k(m(l)-m(l-1))x(l)
\end{equation}
is valid, and vice versa, any non-increasing sequence of the
Boolean variables $ x(1)\ge x(2)\ge\ldots\ge x(k)$  specifies
the label $\mu\in \mcm$ by the formula (\ref{e:mu}).
 By analogy, the feature functions $f_i$ are
represented as  sums
$f_i=\sum_{\tau=1}^L f_i(\tau)$
of non-increasing sequence of the
Boolean variables $ f_i(1)\ge f_i(2)\ge\ldots\ge f_i(L)$.

Let for $l=1\div k$ the vector
$\mbx(l)=(x_1(l), x_2(l),\ldots, x_n(l))$
be Boolean, the vector
$\mbz(l)=(z_1(l),\ldots,z_n(l)$
be with coordinates
\begin{equation*}
z_i(l)=\frac{1}{m(l)-m(l-1)}\sum_{\tau=m(l-1)+1}^{m(l)}f_i(\tau),
\quad (i=1\div n)
\end{equation*}
and the norm
of the vector $|\mbx|=\sum_{i=1}^n|x_i|$, then the
following  proposition is satisfied.
\begin{prop}\label{p:mod}
For any integers   $\nu\in\mcm$ and $f_i\in Z_L$ the equality
\begin{multline*}
|\nu-f_i|=\bigg|m(0)-\sum_{\tau=1}^{m(0)}f_i(\tau)\bigg|+\\
\sum_{l=1}^k\bigg|(m(l)-m(l-1)){\mathbf 1}_{{\big(\ind{\nu\ge m(l)}\big)}}
-\sum_{\tau=m(l-1)+1}^{m(l)}f_i(\tau)\bigg|
+\sum_{\tau=m(k)+1}^{L}f_i(\tau)
\end{multline*}
holds true.
Therefore, for any feature vector $f\in Z_L^V$ and
the Boolean vector
$$\mbx(l)=(1_{\ind{(m_1\ge l)}},1_{\ind{(m_2\ge l)}},
\ldots,1_{\ind{(m_n\ge l)}})$$
the function $U_1(\mbm)$ can be written in the form
\begin{equation}\label{e:u}
U_1(\mbm)=\sum_{l=1}^{k}(m(l)-m(l-1))u(l,\mbx(l)),
\end{equation}
where for $n$-dimensional Boolean vector $\mbb$ functions
\begin{equation*}
u(l,\mbb)=\sum_{i\in V}\la_i|z_i(l)-b_i|+
\sum_{(i,j)\in E}\be_{i,j}
\left|b_i-b_j\right|.
\end{equation*}
\end{prop}

Denote by
\begin{equation}\label{e:lmin}
\check{\mbx}(l)=\arg\min_\mbb u(l,\mbb),\quad (l=1\div k)
\end{equation}
Boolean solutions that minimize the functions
$u(l,\mbb)$.
For two  vectors $v$ and $w$ we will write $v\ge w$  if
all corresponding pairs of their coordinates satisfy the inequality
$v_i\ge w_i,\ (i=1\div n)$ and will write $v\ngeq w$
if there exist at least two different
pairs such that
$v_i\ge w_i$ and $v_j < w_j$.
Note that $\mbz(1)\ge \mbz(2)\ge\ldots\ge\mbz(k)$, and what is more,
for some integer $1\le\varkappa\le k$ their coordinates satisfy  the
following condition

\begin{multline*}
z_i(1)=z_i(2)=\ldots=z_i(\varkappa-1)=1,
0\le z_i(\varkappa)\le 1,\\
z_i(\varkappa+1)=\ldots=z_i(k)=0.
\end{multline*}

It is easy to show that in general the solutions
$\check{\mbx}(l)$ of the (\ref{e:lmin}) are not ordered.
Nevertheless, without fail there is at least one non-increasing
sequence $\check{x}(l)$ of solutions of (\ref{e:lmin}).
\begin{thm}\label{t:mon}
There is a non-increasing sequence
$\check{\mbx}(1)\ge\check{\mbx}(2)\ge\ldots\ge\check{\mbx}(k)$
of solutions of (\ref{e:lmin}).
\end{thm}

\noindent
Some structural properties of the set of
solutions $\check{\mbx}(l)$ are presented in
\begin{cor}\label{c:struct}
For integer $1\le l'<l''\le k$  and
the sequence of vectors $\mbz(1)\ge \mbz(2)\ge\ldots\ge\mbz(k)$
the following properties are valid:

(i) If $\check{\mbx}(l')$ is any solution of (\ref{e:lmin}),
then there is a solution $\check{\mbx}(l'')$ so that
$\check{\mbx}(l')\ge\check{\mbx}(l'')$, and vice versa,
if $\check{\mbx}(l'')$ is any solution of (\ref{e:lmin}),
then there is a solution $\check{\mbx}(l')$ so that
$\check{\mbx}(l')\ge\check{\mbx}(l'')$.

(ii) For each $1\le l\le k$ the set of solutions
$\{\check{\mbx}(l)\}$ has the minimal
$\underline{{\mbx}}(l)$ and the maximal
$\overline{{\mbx}}(l)$ elements.

(iii) The set of minimal and maximal elements are ordered,
i.e. $\underline{{\mbx}}(1)\ge\ldots\ge\underline{{\mbx}}(k)$
and $\overline{{\mbx}}(1)\ge\ldots\ge\overline{{\mbx}}(k)$.
\end{cor}

Sentence \emph{(i)} follows immediately from Theorem \ref{t:mon}.
Sentence \emph{(ii)} is deduced from \emph{(i)} considered for
$l'=l''$, property \emph{(iii)} follows from \emph{(ii)} and definition
of  the minimal and the maximal elements.

For two Boolean vectors $\mbb'$ and $\mbb''$ let the vector
$\uln\mbb=\mbb{'}\wedge\mbb''$, respectively,
$\oln\mbb=\mbb{'}\vee\mbb''$ be with coordinates
$\uln{b}_i=\min\{b_i',b_i''\}$, respectively,
with coordinates $\oln{b}_i=\max\{b_i',b_i''\}$.
If $\check{\mbx}(1),\ldots,\check{\mbx}(k)$
is any unordered sequence of solutions of (\ref{e:lmin})
the ordered sequence of solutions can be derived from it
by the logical operation $\wedge,\vee$ like one-dimensional
variational series
${\mbx}_{(1)}\ge{\mbx}_{(2)}\ge\ldots\ge{\mbx}_{(k)}$.
But it is easy to see the sum of ordered solutions $\check{\mbx}(l)$
is a solution of (\ref{e:m_e}).
\begin{prop}
If $\check{\mbx}(1)\ge\check{\mbx}(2)\ge\ldots\ge\check{\mbx}(k)$
is a sequence of ordered solutions of (\ref{e:lmin}) then
the sum
$$
\wt{\mbm}=m(0)+\sum_{l=1}^k(m(l)-m(l-1))\check{\mbx}(l)
$$
minimizes $U_1(\mbm)$.
\end{prop}

The problem of computing Boolean solutions $\check{\mbx}(l)$
is familiar in the discrete optimization \cite{GPS89}. It is equivalent
to identification of the minimum cuts for  specially built
networks. There are  fast algorithms
to compute them \cite{AMS01,CG94}.

\subsection*{Identification of $\wt{\mbm}_{gaus}$} Now denote the
function to be minimized by
$$
U_2(\mbm)=\sum_{i\in V}\la_i(f_i-m_i)^2+ \sum_{(i,j)\in E}
\be_{i,j}(m_i-m_j)^2, \quad (\mbm\in \mcm^V).
$$
To find a solution $\widehat{\mbm}_{gaus}$ that minimizes
the function $U_2(\mbm)$
the  representation of the vector of cluster  labels
$\mbm$ by the integer valued linear combination of Boolean
vectors is used once more. Then the problem of integer
minimization of the function $U_2(\mbm)$ is reduced
to the problem of Boolean minimization.

Denote for brevity
$ g_i=f_i-m(0),\ (i\in V)$,
$a_l=m(l)-m(l-1),\ (l=1\div k)$.
The vector $\mbm$ can be  represented by the formula (\ref{e:mu})
as the linear combination
$\mbm(0)+\sum_{l=1}^k a_l\mbx(l)$
of Boolean vectors $\mbx(l)=(x_1(l), x_2(l),$ $\ldots, x_n(l))$,
and the function $U_2$ can be written in the form
\begin{equation*}
\sum_{i\in V}\la_i
\bigg(g_i-\sum_{l=1}^k a_l x_i(l)\bigg)^2+
\sum_{(i,j)\in E}\be_{i,j}
\bigg(\sum_{l=1}^k a_l \big(x_i(l)-x_j(l)\big)\bigg)^2.
\end{equation*}
Let $d_{l,\tau}=a_la_\tau,\
(l,\tau=1\div k)$ and $\be_{i,i}=0, (i\in V)$ , then
$
U_2(\mbm)=\sum_{i\in V}\la_i g_i^2+P(\mbx(1),\ldots,\mbx(k)),
$
where the polynomial of Boolean variables $P(\mbx(1),\ldots,\mbx(k))$
after cancellation is written as
\begin{multline*}
P(\mbx(1),\ldots,\mbx(k))=\\
\sum_{i\in V}
\sum_{l=1}^k\bigg[\la_i a_l^2 -2\la_i g_i a_l-\bigg.
\bigg.a_l(m(k)-m(0)-a_l)\sum_{j\in V}(\be_{i,j}+\be_{j,i})
\bigg]x_i(l)+\\
2\sum_{i\in V}\bigg[\la_i+\sum_{j\in V}(\be_{i,j}+\be_{j,i})\bigg]
\sum_{1\le\tau<l\le k}d_{l,\tau}x_i(\tau)x_i(l)+\\
\sum_{(i,j)\in E}\be_{i,j}
\bigg[\sum_{l=1}^k a_l^2 \big(x_i(l)-x_j(l)\big)^2 +
\sum_{l\neq \tau}d_{l,\tau}\big[\big(x_i(l)-x_j(\tau)\big)^2+
(x_j(l)-x_i(\tau))^2\big]\bigg].
\end{multline*}
Note that it has the same points of minimum as $U_2(\mbm)$.
Let us consider another polynomial of Boolean variables
\begin{multline*}
Q(\mbx(1),\ldots,\mbx(k))=\\
\sum_{i\in V}
\sum_{l=1}^k\bigg[\la_i a_l^2 -2\la_i g_i a_l-\bigg.
\bigg.a_l(m(k)-m(0)-a_l)\sum_{j\in V}(\be_{i,j}+\be_{j,i})
\bigg]x_i(l)+\\
2\sum_{i\in V}\bigg[\la_i+\sum_{j\in V}(\be_{i,j}+\be_{j,i})\bigg]
\sum_{1\le\tau<l\le k}d_{l,\tau}x_i(l)+\\
\sum_{(i,j)\in E}\be_{i,j}
\bigg[\sum_{l=1}^k a_l^2 \big(x_i(l)-x_j(l)\big)^2 +
\sum_{l\neq \tau}d_{l,\tau}\big[\big(x_i(l)-x_j(\tau)\big)^2+
(x_j(l)-x_i(\tau))^2\big]\bigg].
\end{multline*}
such that
$Q(\mbx(1),\ldots,\mbx(k))\ge P(\mbx(1),\ldots,\mbx(k))$
and which differs from $P$ by the term
$\sum_{1\le\tau<l\le k}d_{l,\tau}x_i(l)$
in the second line.

Denote by
\begin{equation*}
(\mbq^*(1),\mbq^*(2),\ldots,\mbq^*(k))=
\text{argmin}_{\mbx(1),\mbx(2),\ldots,\mbx(k)}
Q(\mbx(1),\ldots,\mbx(k))
\end{equation*}
any collection of Boolean vectors that minimizes
$Q(\mbx(1),\ldots,\mbx(k))$.
Without fail
$\mbq^*(1)\ge\mbq^*(2)\ge\ldots\ge\mbq^*(L-1)$.
This feature allows expressing solutions of the
initial problem as $\wt\mbm_{gaus}=\sum_{l=1}^{L-1}\mbq^*(l)$.
\begin{thm}\label{t:q=x}
Any collection
$(\mbq^*(1),\mbq^*(2),\ldots,\mbq^*(L-1))$
that minimizes $Q$ forms the nonincreasing sequence.

The polynomials $P$ and $Q$ have the same set of ordered solutions and,
therefore, each solution $\wt\mbm_{gaus}$ is specified
by the formula $\wt\mbm_{gaus}=\sum_{l=1}^{k}\mbq^*(l)$.
\end{thm}

Theorem \ref{t:q=x} allows determination of the
classifier $\wt\mbm_{gaus}$ by the Boolean minimization
of the polynomial $Q$. Unlike $P$ this polynomial can be minimized
directly by the minimum network cut algorithms
\cite{GPS89,AMS01}. The appropriate network is described
in \cite{DAN01}.

\section{Applications}

Here we show several numerical tests as well as a resust of segmentation
of a real 3D US-imade of size $196\times 215\times 301$.

The original 2D gray-scale image in Figure 1a (see next page)
was corrupted by
Gaussian random noise (Figure 1b).
The results of restoration by the extened Ising model
are placed in Figure 1c,1d and by the classical
Ising model are depicted in Figure 1e,1f.

In Figure 2a  the slice of original 3D US-image
of the thyroid gland is depicted.
Its contour that was done by expert is drown in Figure 2a.
The corresponding slice of 3D segmentation of the original image
by the extended Ising model are placed in Figure 2c,d.
The full segmentation of 3D $196\times 215\times 301$ image
takes about 40min of processor Pentium-III 800.

\newpage
\begin{center}
\fbox{\includegraphics[width=5.5cm,height=5.5cm] {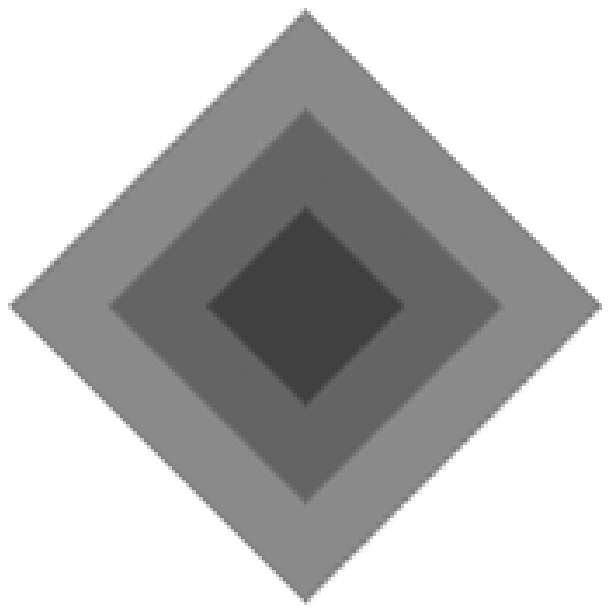}}
\fbox{\includegraphics[width=5.5cm,height=5.5cm]{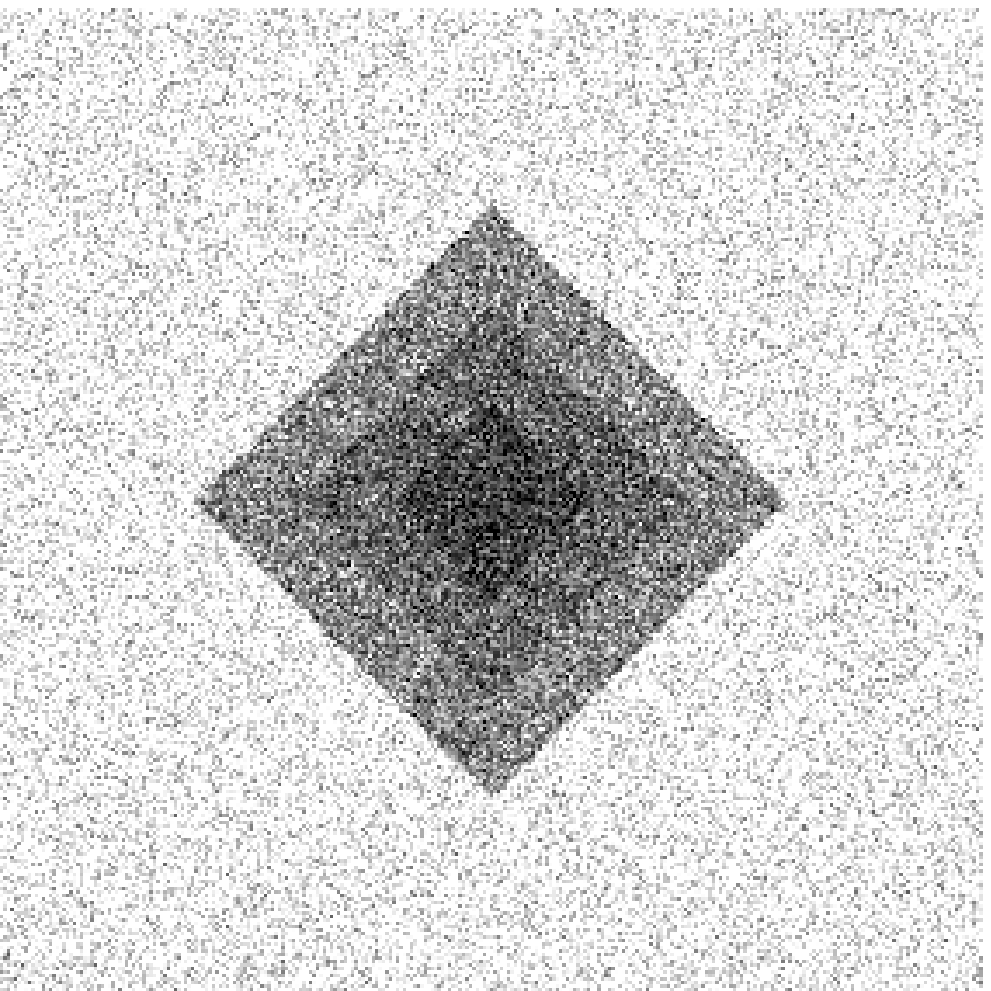}}
\end{center}

\hspace{3cm}(a)\hspace{5.3cm}(b)
\vspace{-0.3cm}

\begin{center}
\fbox{\includegraphics[width=5.5cm,height=5.5cm] {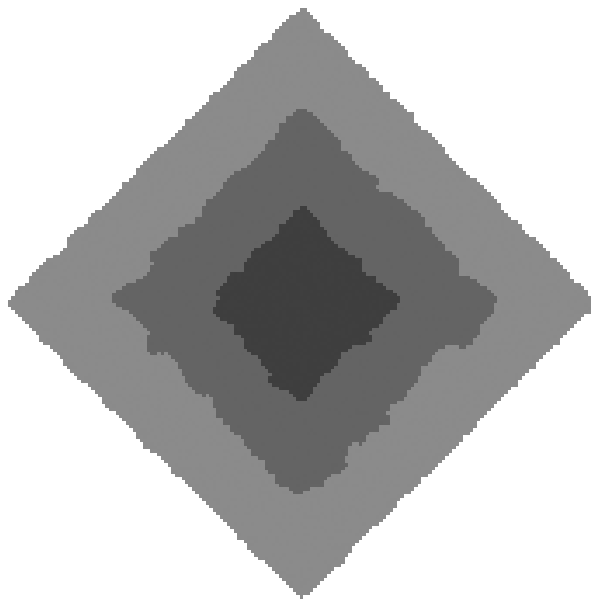}}
\fbox{\includegraphics[width=5.5cm,height=5.5cm]{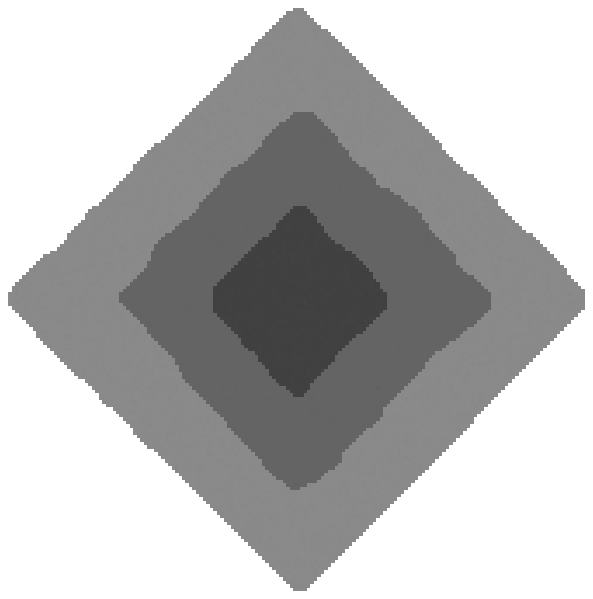}}
\end{center}

\hspace{3cm}(c)\hspace{5.3cm}(d)
\vspace{-0.3cm}

\begin{center}
\fbox{\includegraphics[width=5.5cm,height=5.5cm] {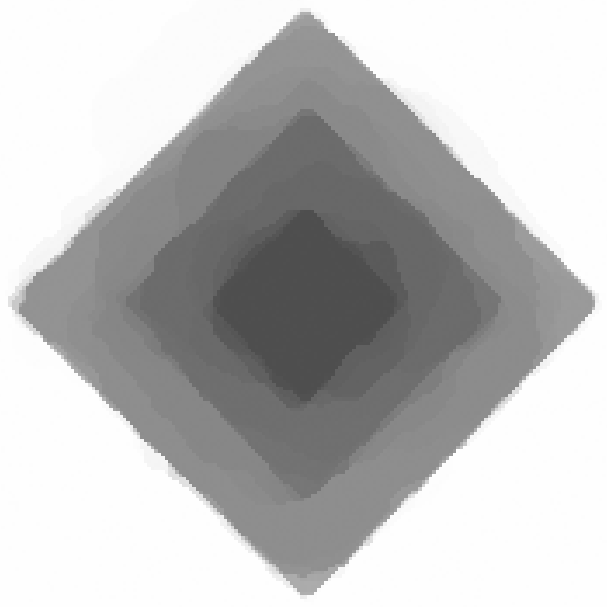}}
\fbox{\includegraphics[width=5.5cm,height=5.5cm]{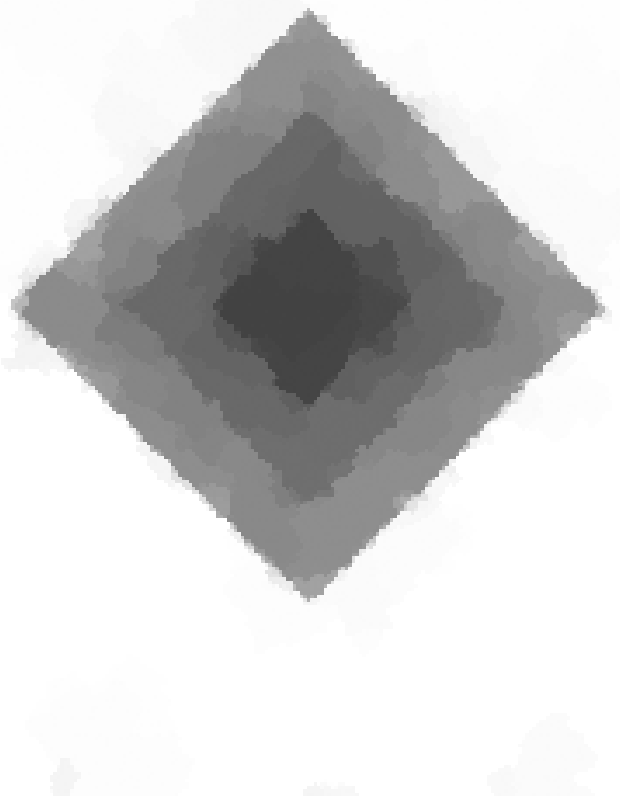}}
\end{center}

\hspace{3cm}(e)\hspace{5.3cm}(f)

\begin{center}
{Figure 1}
\end{center}

\newpage
\begin{center}
\fbox{\includegraphics[width=5.12cm,height=5.69cm] {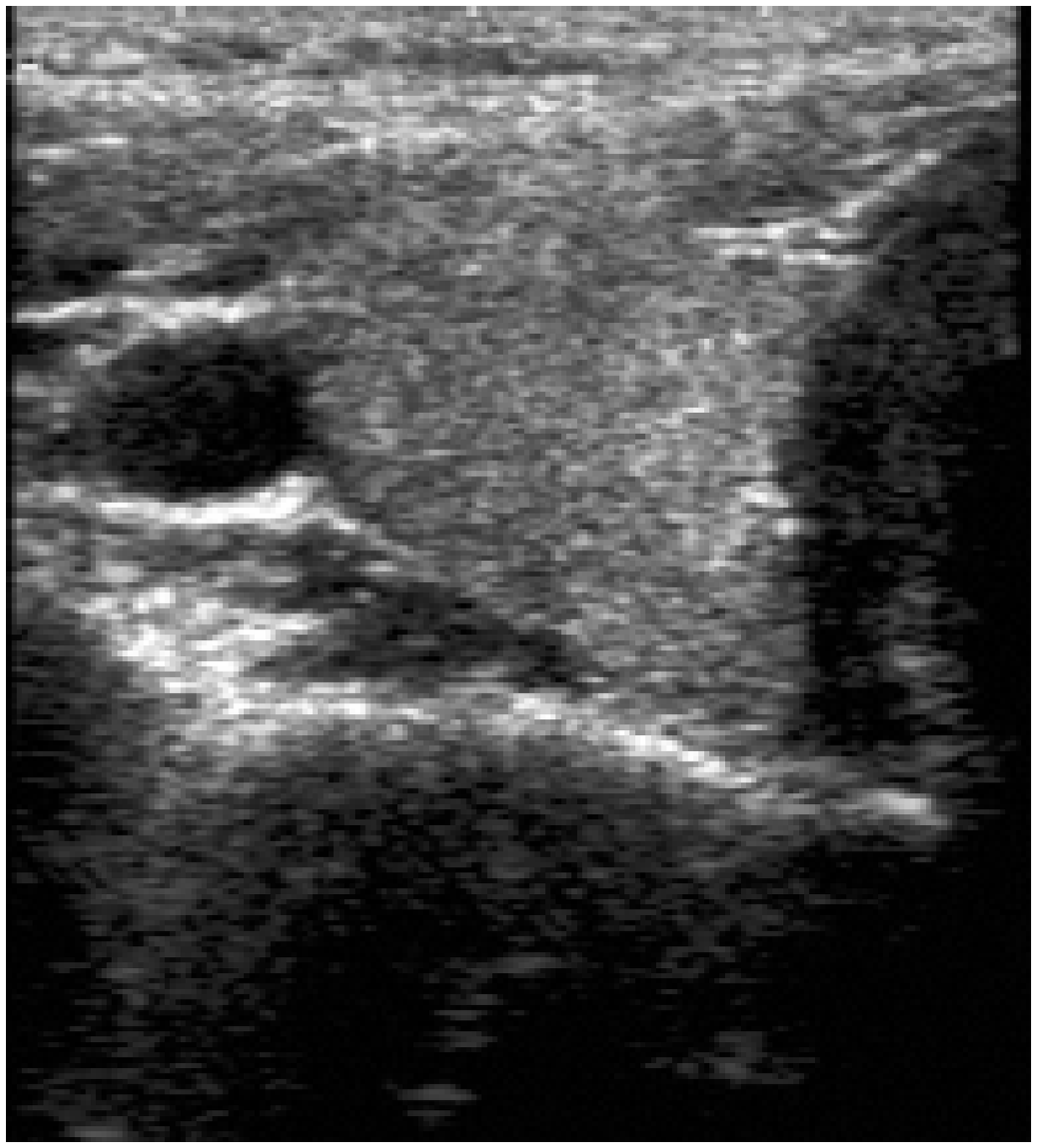}}
\hspace{1cm} \fbox{\includegraphics[width=5.12cm,height=5.69cm]
{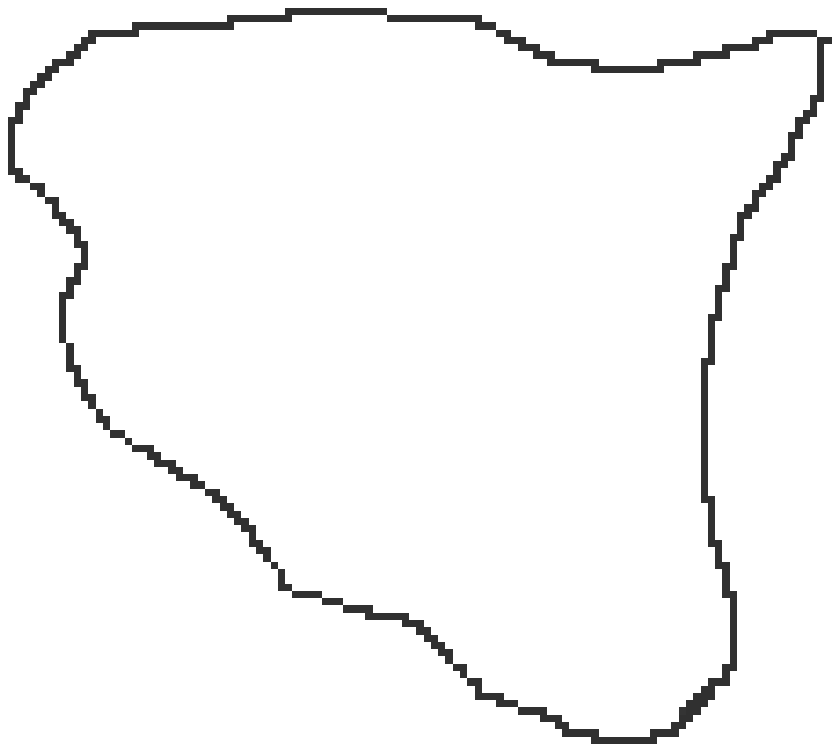}}
\end{center}
\hspace{3cm}(a)\hspace{6.4cm}(b)

\begin{center}
\fbox{\includegraphics[width=5.12cm,height=5.69cm] {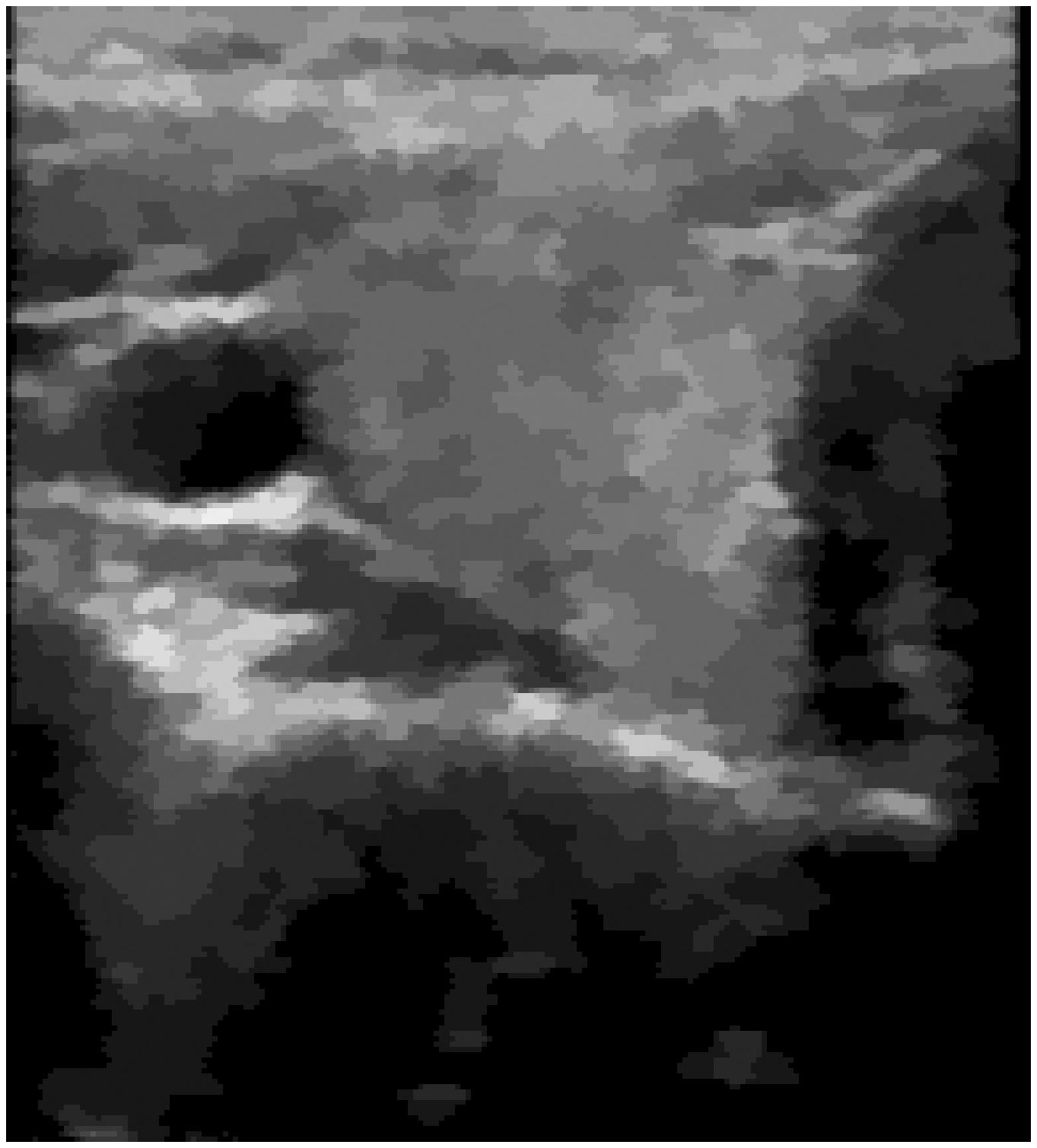}}
\hspace{1cm} \fbox{\includegraphics[width=5.12cm,height=5.69cm]
{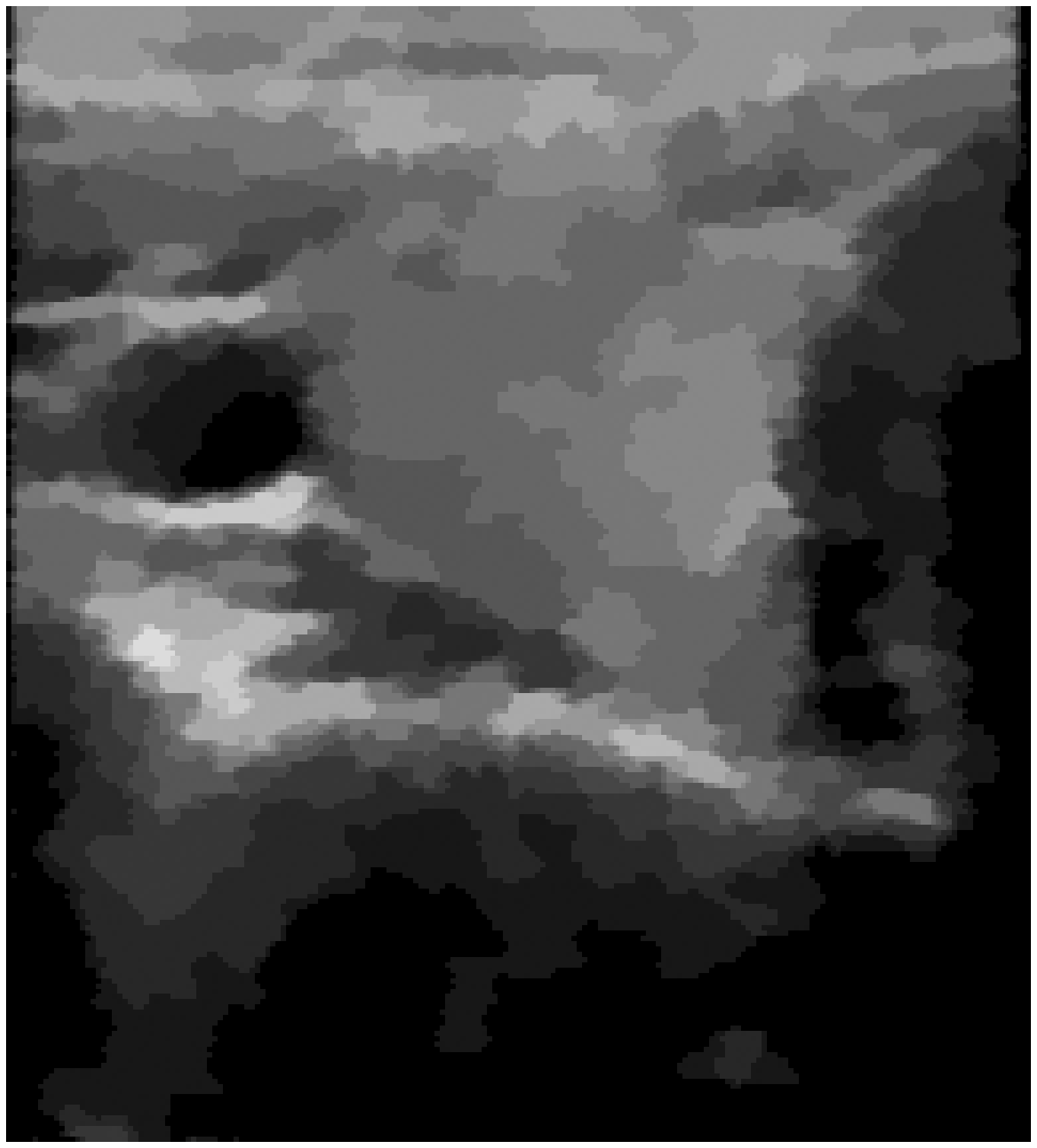}}
\end{center}
\hspace{3cm}(c)\hspace{6.4cm}(d)

\begin{center}
{Figure 2}
\end{center}

\section{Conclusion}

In the paper the  Bayesian methods of segmentation and estimation of
{gray-scale}
and color images are presented. For both of them it is supposed feature
function of images and labels of segments take finite number of rational
(possibly, different) values and they are distributed according to either
the
exponential Gibbs or the Gaussian Gibbs distribution. The numerical
tests showed  the methods developed
allow solution problems of practical segmentation and
Gibbs estimation of images of large sizes.


\begin{thebibliography}{w}

\bibitem{GG84}
S. Geman, D. Geman,
\emph{Stochastic Relaxation, Gibbs Distributions, and the
Bayesian Restoration of Images},
IEEE Trans. on Pattern Anal. and
Mach. Intel. PAMI-6, 6 (1984), 721--741.

\bibitem{Gi95}
B. Gidas,
\emph{Metropolis Type Monte Carlo Simulation Algorithm
and Simulated Annealing}, Topics in Contemp. Probab.
and its Appl., Stochastics Ser., CRC, Boca Raton, Fl,
1995.

\bibitem{Z00}
B.A. Zalesky, \emph{Sthochastic Relaxation for buildind some
classes of piecewise linear regression functions}, Monte-Carlo Meth. and
Apl., 6 (2000), no. 2, 141-157

\bibitem{GPS89}
D.M. Greig , B.T. Porteous, A.H. Seheult,
\emph{Exact Maximum A Posteriori Estimation for Binary
Images}, J. R. Statist. Soc. B.,
58 (1989), 271--279.

\bibitem{BVZ01}
Yu. Boykov, O. Veksler, R. Zabih,
\emph{Fast Approximate Energy Minimization via Graph Cuts},
IEEE Trans. Pattern Anal. and Machine Intel., 23 (2001),
no. 11, 1222--1239.

\bibitem{KZ01}
V. Kolmogorov, R. Zabih,
\emph{What Energy Functions can be Minimized via Graph Cuts?},
Cornell CS Texnical Report, TR2001-1857, 2001.

\bibitem{PRIP01}
B.A. Zalesky,
\emph{Computation
of Gibbs estimates of gray-scale images by discrete
optimization methods}, Proceedings of the Sixth International Conference
PRIP`2001 (Minsk, May 18-20, 2001), 81-85.

\bibitem{DAN01}
B.A. Zalesky,
\emph{ Efficient integer-valued
minimization of quadratic polynomials with the quadratic
monomial $b^2_{i,j}(x_i-x_j)^2$} Dokl. NAN Belarus, 45 (2001),
no. 6, 9-11.


\bibitem{AMS01}
B.A. Zalesky,
\emph{Network Flow Optimization for Restoration of Images},
Preprint AMS, Mathematics ArXiv, 2001, math.OC/0106180.

\bibitem{ABEM89}
S.A. Aivazyan, B.M. Buchshtaber, I.S. Enyukov, L.L. Meshalkin,
\emph{Applied Statistics: Classification and Reducing of Dimension},
Finances and Statistics, Moskow, 1989.

\bibitem{McL92}
G.J. McLachlan,
\emph{Discriminant Analysis and Statistical Pattern Recognition},
John Wiley\&Sons, New York, 1992.

\bibitem{MM99}
V.V. Mottl, I.B. Muchnik,
\emph{Hidden Markov Models in Structure Analysis of Signals},
Fizmatlit, Moskow, 1999.

\bibitem{JH98}
E.E. Juk, Yu.S. Kharin,
\emph{Robustness in Cluster Analysis of Multidimensional data},
Belgosuniversitet, Minsk, 1998,

\bibitem{CG94}
B.V. Cherkassky, A.V. Goldberg,
\emph{On Implementing Push-Relabel Method for the Maximum Flow Problem},
Texnical Report STAN-CS-94-1523, Department of Computer Science,
Stanford University, 1994.

\bibitem{PR75}
J.C. Picard, H.D. Ratliff,
\emph{Minimal cuts and related problems},
Networks, 5 (1975), 357-370.
\end{thebibliography}
\end{document}